\documentclass[12pt]{amsart}

\usepackage{amssymb,amsthm,amsfonts,graphicx,color,comment,amsmath}
\usepackage{enumerate}

\newtheorem{thm}{Theorem}[section]
\newtheorem*{thm*}{Theorem}
\newtheorem{lemma}[thm]{Lemma}

\newtheorem{prop}[thm]{Proposition}

\theoremstyle{definition}
\newtheorem{definition}[thm]{Definition}
\newtheorem{remark}[thm]{Remark}

\newtheorem{example}[thm]{Example}

\def\R{{\mathbb R}}
\def\Z{{\mathbb Z}}
\def\N{{\mathbb N}}
\def\C{{CAT(0) }}

\def\E{{\mathbb E}}

\title{Limit groups are CAT(0)}
\author{Emina Alibegovi\'c and Mladen Bestvina}
\address{Department of Mathematics, University of Michigan, Ann Arbor,
  Michigan} 
\email{eminaa@umich.edu}
\address{Department of Mathematics, University of Utah, Salt Lake
  City, Utah}
\email{bestvina@math.utah.edu}
\begin{document}

\begin{abstract}
We prove that every limit group acts geometrically on a CAT(0) space
with the isolated flats property.
\end{abstract}
\maketitle
\section{Introduction}

A group is said to be a CAT(0) group if it acts geometrically,
i.e. properly discontinuously and cocompactly by isometries, on a
CAT(0) space. One should think of a \C space as a geodesic metric
space in which every geodesic triangle is at least as thin as its
comparison triangle in Euclidean plane. For basic facts about \C spaces
and groups a general reference is \cite{MR2000k:53038}. We will be
interested mainly in geodesic spaces which are locally \C spaces
(non-positively curved spaces), i.e. every point has a neighborhood
which is a \C space. All of our (locally) CAT(0) spaces will be
proper. In this paper we show that the class of groups known as {\it
limit groups} is CAT(0), thus answering a question circulated by
Z. Sela.

Limit groups arise naturally in the study of equations over free
groups. Historically they appear under different names: $\exists$-free
groups \cite{MR91f:03077}, fully residually free groups
\cite{MR2000b:20032a}, \cite{MR2000b:20032b},
$\omega$-residually free groups \cite{zlil1}. The term ``limit group''
first appears in \cite{zlil1} and reflects the topological approach to
the subject in which these groups arise in the limit of a sequence of
homomorphisms $G\to F$ ($F$ a free group) interpreted as a sequence of
$G$-trees. Thus limit groups act isometrically on $\R$-trees (see also
\cite{MR93j:20059}) and the Rips machine can be brought to bear. For
more information on limit groups the reader is referred to the 
cited papers above and to the expositions
\cite{champetier-guirardel},\cite{paulin},\cite{bf:notes}.

The quickest definition, in the incarnation of $\omega$-residually
free groups, is that a finitely generated group $L$ is a limit group
provided that for every finite subset $W\subset L$ there is a
homomorphism $L\to F$ to a free group which is injective on $W$. A
more constructive definition will be given in the Appendix (see
Proposition \ref{CLG}). For the
convenience of the reader we summarize the basic 
properties of limit groups in the following
proposition. Properties (1) and (3)-(5) follow easily from the
definition, (2) and (6)-(9) can be found in the above papers (and can
be proved quickly from the constructive definition), and (10)
was proved independently by F. Dahmani \cite{dahmani} and
E. Alibegovi\' c \cite{emina}.

\begin{prop}\label{basic}
Let $L$ be a limit group.
\begin{enumerate}
\item[(1)] $L$ is torsion-free.
\item[(2)] Every abelian subgroup of $L$ is finitely generated.
\item[(3)] Every nontrivial abelian subgroup is contained in a unique
  maximal abelian subgroup.
\item[(4)] If two nontrivial commuting elements are conjugate then
  they are equal.
\item[(5)] Roots are unique, i.e. $x^n=y^n$ with $n\neq 0$ implies
  $x=y$.
\item[(6)] $L$ is finitely presented, and in fact has a finite
  Eilenberg-MacLane space. 
\item[(7)] $L$ is coherent, i.e. every finitely generated subgroup is
  finitely presented (and it is a limit group).
\item[(8)] If $L$ is not free then its cohomological dimension is
  $\max(2,n)$ where $n$ is the maximal rank of an abelian subgroup. 
\item[(9)] There are finitely many conjugacy classes of maximal
  abelian subgroups of rank $\geq 2$.
\item[(10)] $L$ is relatively hyperbolic with respect to the collection
  of maximal abelian subgroups of rank $\geq 2$.
\end{enumerate}
\end{prop}

\begin{thm*} Limit groups are CAT(0).
\end{thm*}

Here is an outline of the argument, contained in Sections 2 and 3.
An important class of limit groups consists of {\it
$\omega$-residually free towers} \cite[Definition 6.1]{zlil1} (see
also \cite{MR2000b:20032b}). These are groups that have explicit and
inductive description in terms of gluing building blocks (graphs,
surfaces, tori) in a certain way. We recall this definition in Section
\ref{towers}. It is easy to see that $\omega$-residually free towers
are CAT(0). We then make use of the fact
\cite{zlil2},\cite{MR2000b:20032b}, whose proof due to Sela we outline
in the Appendix, that every limit group occurs as a
finitely generated subgroup $L$ of an $\omega$-residually free tower
$G$. The associated \C space is then constructed as a
``core'' of the covering space of $G$ associated to $L$.

A closer analysis of the construction shows that the 
CAT(0) space we built satisfies the {\it isolated flats property}
\cite{kruska}. This analysis is contained in Section 4. After this
paper was written, the results of Hruska-Kleiner \cite{hruska-kleiner}
became available. In particular, it follows from \cite[Theorem
  1.2.1]{hruska-kleiner} and Proposition \ref{basic}(10) that {\it
  every} CAT(0) space on which a limit group acts geometrically must
have isolated flats.

\nocite{drutusapir}

We thank Noel Brady for pointing out the simplification of the
statement of Main Theorem. 

\section{Geometrically coherent locally CAT(0) complexes}

Recall that a group is {\it coherent} if its every finitely generated subgroup
is finitely presented. We now introduce a \C version of
coherence. 

\begin{definition}
Let $X$ be a connected locally \C space. A subspace $C$ of $X$ is a
{\it core} of $X$ if it is compact, locally \C (with respect to the
induced path metric), and inclusion $C\to X$ is a homotopy
equivalence.
\end{definition}

Note that if $C$ is connected and locally \C then it is
aspherical. Thus inclusion $C\to X$ is a homotopy equivalence if and
only if it induces an isomorphism in $\pi_1$.

Recall that a celebrated theorem of Peter Scott \cite{MR0326737}
asserts that any connected 3-manifold $M$ with finitely generated
$\pi_1(M)$ contains a compact
submanifold $N$ such that $\pi_1(N)\to\pi_1(M)$ is an isomorphism, and
in particular $\pi_1(M)$ is finitely presented. The compact
submanifold $N$ is known as a {\it Scott core} of $M$. Scott's theorem is a
motivation for our approach to proving that limit groups are $CAT(0)$.

\begin{definition}
A connected locally \C space $Y$ is {\it geometrically coherent} if for every
covering space $X\to Y$ with $X$ connected and $\pi_1(X)$ finitely
generated and every compact subset $K\subset X$ it follows that $X$
contains a core $C\supset K$.
\end{definition}

D. Wise has studied this property for 2-complexes in the presence of
sectional curvature conditions \cite{wise:seccurv}. 

\begin{example}
Any compact hyperbolic surface with totally geodesic boundary is
geometrically coherent. Moreover, a core can always be chosen to be
convex (see e.g. \cite[4.73,4.75]{misha}).
\end{example}

\begin{example}
A flat torus $T=\R^n/\Lambda$ is geometrically
coherent. Indeed, any connected covering space admits a metric
splitting $T'\times E$ where $T'$ is a flat torus or a point and $E$
is Euclidean space $\R^k$ or a point. Thus $T'\times\{point\}$ is a
(convex) core.
\end{example}

\begin{example}
Without the curvature assumptions there are examples due to Wise
\cite{MR1934009} of finite 2-complexes $X$ and covering spaces $\tilde
X\to X$ such that $\pi_1(\tilde X)$ is finitely presented and yet
$\tilde X$ does not contain a compact subset $C$ with
$\pi_1(C)\to\pi_1(\tilde X)$ an isomorphism. We also point out the the
product of two finite graphs with nonabelian fundamental groups is a
2-dimensional locally $CAT(0)$ complex which is not coherent, much
less geometrically coherent. It appears to be unknown whether a
locally $CAT(0)$ complex which is coherent is necessarily
geometrically coherent. 
\end{example}

Many examples of geometrically coherent spaces can be constructed by
repeatedly applying the following gluing theorem. The annulus $S^1\times
[0,1]$ is endowed with a flat metric where $S^1$ has suitably chosen
total length.

\begin{thm}\label{main}
Let $M$ and $N$ be geometrically coherent locally \C spaces. Let $Y$ be
the space obtained from the disjoint union of $M$, $N$ and a finite
collection of annuli $S_t^1\times \{0\}$
$$M\sqcup \sqcup_t S_t^1\times [0,1]\sqcup N$$ by gluing
$S_t^1\times\{0\}$ to a local geodesic in $M$ by a local isometry and
gluing $S_t^1\times\{1\}$ to a local geodesic in $N$ by a local
isometry.

Then $Y$ is geometrically coherent.
\end{thm}

In the proof we shall need the following fact.

\begin{prop}(\cite[Gluing with a tube, 11.13]{MR2000k:53038})\label{gluing}
Let $X$ and $A$ be locally CAT(0) metric spaces. If $A$ is compact and
$\phi, \psi: A \to X$ are locally isometric immersions, then the
quotient of $X \sqcup (A\times [0,1])$ by the equivalence relation
generated by $(a,0) \sim \phi(a); (a,1) \sim \psi(a), \forall a \in A$
is locally CAT(0).
\end{prop}

\begin{proof}[Proof of Theorem \ref{main}]
Let $p:X \to Y$ be a covering space with $X$ connected and $\pi_1(X)$
finitely generated and let $K \subset X$ be a given compact set. Note
that $X$ is a graph of spaces with vertex spaces the components of
$p^{-1}(M)$ and $p^{-1}(N)$ and with edge spaces the components of
$p^{-1}(S_t^1\times \{\frac 12\})$. Let $\Gamma$ be the graph of groups
associated to $X$, so that $\pi_1(X)=\pi_1(\Gamma)$. Thus $\Gamma$ is
a bipartite graph with two types of vertices which we will call
$M$-type and $N$-type. All edge groups are either trivial or infinite
cyclic. The graph $\Gamma$ may be infinite, so consider an exhaustion
of $\Gamma$ by a chain of connected finite subgraphs. This induces an
exhaustion of $\pi_1(X)$ by a chain of subgroups. Since $\pi_1(X)$ is
finitely generated this chain stabilizes, so there is a finite
subgraph $\Gamma_0\subset \Gamma$ with
$\pi_1(\Gamma_0)=\pi_1(\Gamma)$. Let $X_0$ be the subspace of $X$
corresponding to $\Gamma_0$. $X_0$ may not contain the given compact
set $K$, but $K$ will certainly be contained in a finite subgraph of
spaces of $X$, so we will simply enlarge $\Gamma_0$ and $X_0$ so as to
contain $K$. Each edge of $\Gamma_0$ corresponds to an annulus or a
strip ($\R \times [0,1]$). One boundary component of each annulus or
strip is identified with a circle or a line in an $M$-type vertex and
the other is identified with a circle or a line in an $N$-type vertex.
As a graph of groups, $\Gamma_0$ is finite with trivial or cyclic edge
groups, and $\pi_1(\Gamma_0)=\pi_1(X)$ is finitely generated. It
follows that all vertex groups are finitely generated so we are in
position to apply our assumptions about the existence of cores in
these vertex spaces.

We will now build a core $C$ in $X_0$ (and this will be a core in $X$
as well). The idea is to take the union of cores in vertex spaces,
connecting annuli, plus rectangles in connecting strips. The cores and
the rectangles have to be chosen with care so that the union contains $K$
and so that the intersections between rectangles and neighboring cores
is contractible.

Recall that $X_0$ contains finitely many vertex and edge spaces.
Start by selecting a rectangle $R(S)=I(S)\times [0,1]$ inside
each strip $S=\R\times [0,1]\subset X_0$ so that
\begin{enumerate}
\item[(1)] $R(S)\supset S\cap K$.
\end{enumerate}

Now for each vertex space $V\subset X_0$ choose a core $C(V)$ such
that

\begin{enumerate}
\item[(2)] $C(V)\supset V\cap K$,
\item[(3)] $C(V)\supset V\cap R(S)$ for every strip $S\subset X_0$, and
\item[(4)] $C(V)\supset V\cap A$ for every annulus $A\subset X_0$.
\end{enumerate}

Let $C$ be the union of the following spaces:
\begin{itemize}
\item $C(V)$ for every vertex space $V\subset X_0$,
\item $A$ for every annulus $A\subset X_0$, and
\item $R(S)$ for every strip $S\subset X_0$.
\end{itemize}

Then $C$ is a compact space and it has the same
fundamental group as $X_0$. That $C$ is locally CAT(0) follows 
from Proposition \ref{gluing}.
\end{proof}

\begin{remark}
In the statement of the Theorem one can replace annuli $S^1\times
[0,1]$ by products $T\times [0,1]$ of flat tori $T$. The proof is similar. 
\end{remark}

\section{$\omega$-residually free towers}\label{towers}

\begin{definition}\cite[Definition 6.1]{zlil1}\label{AQ}
A {\it height 0 $\omega$-rft} is the wedge of finitely
many circles, tori and closed hyperbolic surfaces excluding the
surface of Euler characteristic $-1$.

Assume that inductively the notion of a height $n-1$
$\omega$-rft has been defined. A height $n$
$\omega$-rft $Y_n$ is obtained from a height $n-1$
$\omega$-rft $Y_{n-1}$ by gluing a building
block. 
\begin{itemize}
\item[(A)] (Abelian block) $Y_n=Y_{n-1}\sqcup S^1\times [0,1]\sqcup T^m/\sim$
  where $T^m$ is the $m$-torus, $S^1\times \{1\}$ is identified with a
  coordinate circle in $T^m$ and $S^1\times\{0\}$ is identified with a
  nontrivial 
  loop in $Y_{n-1}$ that generates a maximal abelian subgroup in
  $\pi_1(Y_{n-1})$. 
\item[(Q)] (Quadratic block) $Y_n=Y_{n-1}\sqcup \Sigma \sqcup_t
  S_t^1\times [0,1]/\sim$ where $\Sigma$ is a connected compact
  hyperbolic surface with totally geodesic boundary and $\chi\leq -2$
  or a punctured torus\footnote{these are precisely the hyperbolic
    surfaces that support pseudoAnosov homeomorphisms or equivalently
    two intersecting two-sided simple closed curves}, $S_t^1\times
  \{0\}$ is identified with a boundary component of $\Sigma$, and
  $S_t^1\times \{1\}$ is identified with a homotopically nontrivial
  loop in $Y_{n-1}$. A further requirement is that there exists a
  retraction
$$r:Y_n\to Y_{n-1}$$ such that the restriction $r:\Sigma\to Y_{n-1}$
  has nonabelian image in $\pi_1$.
\end{itemize}
A group $\Gamma$ is called an $\omega$-rft if it is the fundamental
group $\Gamma=\pi_1(Y)$ of a
space which is an $\omega$-rft.
\end{definition}

We took certain liberties stating this definition. In particular, Sela
allows $Y_n$ to be obtained from $Y_{n-1}$ by attaching more than one
building block and also allows wedging tori and closed
surfaces. However, it is easy to see that both definitions define the
same class of groups.

\begin{lemma}
Every $\omega$-rft is geometrically coherent.
\end{lemma}

\begin{proof}
Every $\omega$-rft is given a locally \C metric by induction on
height. When a surface is glued it is endowed with a hyperbolic metric with
totally geodesic boundary. We use the fact that a hyperbolic
structure can be given with prechosen lengths of boundary components
\cite[expos\'e 3, \S II]{flp},
so that these lengths match the lengths of circles to which they are
glued. The statement now follows from Theorem \ref{main} and induction
on height.
\end{proof}

The importance of $\omega$-rft's is that any limit
group embeds in one. The main result of \cite{MR2000b:20032b} is an
essentially equivalent statement but expressed in a different
language. We outline Sela's construction in the Appendix.

\begin{thm}\cite[1.11,1.12]{zlil2}\label{embedding}
Every limit group is isomorphic to a finitely generated subgroup of an
$\omega$-residually free tower.
\end{thm}

This in turn immediately implies the result we are aiming for.

\begin{thm}\label{cat0}
Limit groups are CAT(0).
\end{thm}

\section{Isolated flats}

In this section we give a combination theorem for CAT(0) spaces with
isolated flats which we will then use to show that limit groups act 
geometrically on such spaces. 
We recall the definition of isolated flats property given by C. Hruska in
\cite{kruska}. A flat in a CAT(0) space $X$ is an
isometric embedding of a Euclidean space $\E^k$ into $X$ for some $k
\geq 2$. A half-flat is an isometric embedding of $\E^{k-1}\times
[0,\infty)$.

\begin{definition}\cite[3.2]{kruska}\label{iflats}
A CAT(0) space $X$ has isolated flats property if it contains a family
$\mathcal F$ of flats so that the following are satisfied: 

\begin{enumerate}
\item There is a constant $B$ so that every flat in $X$ is contained
  in the Hausdorff $B$-neighborhood of some flat $F \in \mathcal{F}$. 
\item There exists $\phi: \R_+ \to \R_+$ so that for every pair of
  distinct flats $F_1, F_2 \in \mathcal F$ and for every $k \geq 0$
  the intersection of Hausdorff $k$-neighborhoods $N_k(F_1)\cap
  N_k(F_2)$ has diameter at most $\phi(k)$.  

\end{enumerate}
\end{definition}

It is not hard to show that the family $\mathcal F$ can be
chosen so that it is invariant under all isometries of the space
$X$.

Suppose a group $G$ acts geometrically on a metric space $X$. A
$k$-flat $F$ in $X$ is {\it periodic} if there is a free abelian
subgroup $A<G$ of rank $k$ that acts on $F$ by translations with a
quotient a $k$-torus.  We will also say that a line (i.e. a biinfinite
geodesic) $\ell$ in $X$ is periodic if there is an element $g
\in G$ which acts on it as a nontrivial translation. 
If $G$ acts geometrically on a CAT(0) space
$(X,\mathcal F, \phi)$ with isolated flats, then every flat $F\in
\mathcal F$ is periodic (\cite[3.7]{kruska}).

We say that a line $\ell$ in $X$ is parallel to a flat or a line $F
\subset X$ if $\ell$ is contained in a Hausdorff neighborhood of $F$.
If a line is parallel to another line, then the two cobound a
strip (\cite[II.2.13]{MR2000k:53038}). If a line $\ell$ is parallel to
a flat $F$ then $F$ contains a line parallel to $\ell$.

\begin{lemma}\label{psi}
Let $(X,\mathcal F, \phi)$ be a CAT(0) space with isolated
flats on which a group $G$ acts geometrically. Let
$\ell$ be a periodic geodesic in $X$ not parallel to any flat $F \in
\mathcal F$. There exists a function $\psi:\R_+ \to \R_+$ so that for
every $k \geq 0$ and every $F\in \mathcal F$ the diameter of
$N_k(\ell) \cap N_k(F)$ is no bigger than $\psi(k)$.
\end{lemma}

\begin{proof}
We will assume without loss of generality that $\mathcal F$ is
$G$-invariant. 
Suppose such a function does not exist. Then there exist $k \in \N$
and a sequence of flats $F_i \in \mathcal F$ so that diam$(N_k(\ell)
\cap N_k(F_i)) > i$. Let $\ell_i=N_{2k}(F_i)\,\cap\,\ell$. Note that
as $i \to \infty$ the length $l_i$ of $\ell_i$ (possibly infinite)
also tends to infinity. Fix a point $p\in\ell$. Since $\mathcal F$ is
preserved by $G$ we may, if necessary, replace each $F_i$ by a
translate $g_i(F_i)$ for some element $g_i\in G$ that acts as a
translation on $\ell$ so that for large $i$ we have $p\in\ell_i$ and,
moreover, the distance between $p$ and the
endpoints (if any) of $\ell_i$ goes to infinity. But then
$N_{2k}(F_i)\cap N_{2k}(F_j)$ contains a large neighborhood in $\ell$
of $p$ (for large $i,j$) which implies that $F_i=F_j$. Thus the
sequence $F_i$ eventually consists of a single flat and $\ell$ is
parallel to it. 
\end{proof}

We are interested in finite graphs of spaces that are locally $CAT(0)$
  and whose edge spaces are geodesic circles or segments.
 To reduce
  repetition we will call such a space $X$ a {\it special CAT(0) graph
  of spaces}. Of course, the universal cover $\tilde X$ of $X$ has an
  induced decomposition as a tree of spaces whose edge spaces are
  lines or segments. In what follows we adopt the following
  terminology. Let $\pi:\tilde X\to T$ be the natural map to the tree
  that collapses each vertex space $\tilde X_v$ to the associated
  vertex $v$ and collapses products $\tilde X_e\times [0,1]$ to edges
  of $T$. If $\tilde X_v$ is a vertex space, we define the {\it
  extended vertex space} $\tilde X_v^+$ to be the inverse image of the
  star of $v$ with respect to the barycentric subdivision of $T$. Thus
  $\tilde X_v^+$ is the union of $\tilde X_v$ with products $\tilde
  X_e\times [0,\frac 12]$ (or $\tilde X_e\times [\frac 12,1]$) for
  each $\tilde X_e\times [0,1]$ attached to $\tilde X_v$. An edge
  space $\tilde X_e$ is identified with $\tilde X_e\times\{\frac
  12\}\subset \tilde X_e\times [0,1]$. Thus an edge space which is
  incident to a vertex space $\tilde X_v$ is contained in the
  associated extended vertex space $\tilde X_v^+$.

\begin{thm}\label{combination}
Let $X$ be a special CAT(0) graph of spaces whose vertex spaces come
in two colors: G (for ``good'') and B (for ``bad''). We consider the
induced coloring of the vertex spaces of $\tilde X$. Assume:
\begin{enumerate}
\item[(0)] each edge space $\tilde X_e$ isometric to a line is
  incident to at least one G vertex,
\item[(1)] distinct edge spaces in an extended G vertex space are not
  parallel to each other, and
\item[(2)] an edge space in an extended G vertex space $\tilde X_v^+$ does not
bound a half-flat in $\tilde X_v^+$.
\end{enumerate}
If the universal cover of each vertex space $X_v$ has the isolated
flats property, then so does the universal cover $\tilde X$ of $X$.
\end{thm}

Note that (2) implies
\begin{enumerate} 
\item[(3)] an edge space in an extended G vertex space $\tilde
  X_v^+$ is not parallel to a flat in $\tilde
  X_v^+$.
\end{enumerate}

In the proof we will need two lemmas.

\begin{lemma}\label{half-flats0}
Let $Z$ be a CAT(0) space, $u,v$ two isometries of $Z$ such that
$<u,v>\subset Isom(Z)$ is discrete, and $A_u,A_v$ are axes of $u,v$
respectively. If for some $k>0$ the intersection $N_k(A_u)\cap
N_k(A_v)$ has infinite diameter then
some nontrivial powers of $u$ and $v$ coincide, and in particular $A_u$ and
$A_v$ are parallel.
\end{lemma}

\begin{proof}
Let $x_n\in A_u$, $y_n\in A_v$ be sequences of points going to
infinity with $d(x_n,y_n)\leq 2k$, $n=0,1,2,\cdots$. Denote by
$|u|,|v|$ the translation lengths of $u,v$ respectively. For each $n$
there are powers $u^{\alpha_n},v^{\beta_n}$ of $u,v$ taking $x_0,y_0$
to within $|u|,|v|$ of $x_n,y_n$. Thus $v^{-\beta_n}u^{\alpha_n}$
takes $x_0$ to within $2k+|u|+|v|$ of $y_0$. By the discreteness
assumption (along with the blanket assumption that $Z$ is proper)
there are only finitely many elements of $<u,v>$ with this property
and we conclude $$v^{-\beta_n}u^{\alpha_n}=v^{-\beta_m}u^{\alpha_m}$$
for many $\beta_m\neq\beta_n$ and the claim follows.
\end{proof}

\begin{lemma}\label{psi'}
Under the hypotheses of Theorem \ref{combination}
there is a function $\psi':\R_+\to\R_+$ so that for every extended G
  vertex space $\tilde X_v^+$ and any two distinct edge spaces $\tilde
  X_e,\tilde X_f\subset \tilde X_v^+$
  the diameter of $N_k(\tilde
  X_e)\cap N_k(\tilde
  X_f)$ is no bigger than $\psi'(k)$.
\end{lemma}

 \begin{proof} Suppose the statement fails. Then there is $k>0$ and a
  sequence of pairs of edge spaces in extended G vertex spaces whose
  $k$-neighborhoods have larger and larger intersections.  Since there
  are only finitely many orbits of extended vertex spaces, and in each
  extended vertex space only finitely many orbits of edge spaces, we
  may assume that one of the edge spaces in these pairs is a fixed
  edge space $\tilde X_e$.  Now we can translate the other edge spaces
  in these pairs by isometries stabilizing $\tilde X_e$ (as in the
  proof of Lemma \ref{psi} except that flats are replaced by other
  edge spaces). Since a compact set intersects only finitely many edge
  spaces, we conclude that there is another edge space $\tilde X_f$ so
  that $N_k(\tilde X_e)\cap N_k(\tilde X_f)$ has infinite
  diameter. Now we have a contradiction to Lemma \ref{half-flats0}.
\end{proof}

\begin{proof}[Proof of Theorem \ref{combination}]

Let $q: \tilde X \to X$ be the universal covering space. We need to
find a family of flats $\mathcal F$ in $\tilde X$ and the function
$\phi: \R_+ \to \R_+$ so that for any two flats $F_1, F_2 \in \mathcal
F$ we have diam$(N_k(F_1)\cap N_k(F_2)) \leq \phi(k)$.

A fixed connected component of $q^{-1}(X_v)$ is the universal cover of
$X_v$ and, by our assumption, has a family of flats $\mathcal F_v$ and
a function $\phi_v$ as in Definition \ref{iflats}. The family of flats
for any other connected component of $q^{-1}(X_v)$ is just a translate
of $\mathcal F_v$ by an appropriate element of $\pi_1(X)$. It is
therefore immediate that the function we seek in Definition
\ref{iflats}(2) for these families will be $k\mapsto\phi_v(2k)+2k$ (a
point $x$ of $\tilde X$ within $k$ of both $F_1\subset \tilde X_v$ and
$F_2\subset \tilde X_v$ is within $k$ of a point in $\tilde X_v$ which
is within $2k$ of both $F_1$ and $F_2$ -- this point can be taken to
be the first point of $\tilde X_v$ on the shortest geodesic from $x$
to $F_1$).

Let $\mathcal F=\{gF: g \in \pi_1(X), F \in \mathcal F_v, \text{for
  some vertex} \ v \in \Gamma \}$. 
We will show that this family satisfies
  the requirements from the definition of isolated flats property for
  the space $\tilde X$. 

{\bf Claim.} Every flat in $\tilde X$
  is contained in some vertex space. 

Indeed, suppose there is a flat $F \subset \tilde X$, necessarily
  2-dimensional, which is not contained in any vertex space. Consider
  the intersection of $F$ with the collection of edge spaces in
  $\tilde X$. This intersection consists of parallel lines in $F$. Let
  $\ell$ be one of these lines. Thus $\ell$ is an edge space $\tilde
  X_e$ and one of the two complementary regions in $F$ adjacent to
  $\ell$ is a strip or a half-plane in an extended G vertex space,
  violating (1) or (2). The claim is proved.

  Let $F_1$ and $F_2$ be two flats in $\mathcal F$ contained in two
  distinct vertex spaces $\tilde X_v$ and $\tilde X_w$ respectively.
  There is an edge path $\tilde X_{e_1}, \ldots, \tilde X_{e_n}$ in
  $\tilde X$ between $\tilde X_v$ and $\tilde X_w$. If $\tilde
  X_{e_i}$ is a segment, for some $i$, then diam$(N_k(F_1)\cap
  N_k(F_2)) \leq $ diam$(N_k(\tilde X_{e_i}))\leq\mbox{diam}(\tilde
  X_{e_i})+2k$. We now assume that all edge spaces in this edge path
  are lines. According to the assumption (0) there are two cases.
  
  {\bf Case 1.}  One of $\tilde X_v$ or $\tilde X_w$, say $\tilde
  X_v$, is a G vertex. Let $e=e_1$. By Lemma \ref{psi} (which applies
  by (3)), there is a
  function $\psi_e$ so that the diameter of $Q=N_{2k}(F_1)\cap 
  N_{2k}(\tilde X_e \times \{0\}) \cap \tilde X_v$ is $ \leq
  \psi_e(2k)$.  Note that $N_k(F_1)\cap N_k(F_2)$ is contained in the
  $k$-neighborhood of $Q$ (if $x\in N_k(F_1)\cap N_k(F_2)$ is on the
  $\tilde X_v$-side of $\tilde X_e \times \{0\}$ then the first point
  of intersection with $\tilde X_e \times \{0\}$ of the shortest
  geodesic joining $x$ to $F_2$ is in $Q$, and similarly for the other
  side), hence diam$(N_k(F_1)\cap
  N_k(F_2)) \leq \psi_e(2k)+2k$.
  
  {\bf Case 2.} One of the interior vertex spaces of the edge path is
  a G vertex. Say this is the vertex space $\tilde X_u$ and the
  adjacent edge spaces are $\tilde X_{e_i}$ and $\tilde X_{e_{i+1}}$,
  so that $\tilde X_{e_i}\subset \tilde X_u^+$ and $\tilde
  X_{e_{i+1}}\subset \tilde X_u^+$. By Lemma \ref{psi'} we have that
  $Q=N_{2k}( \tilde X_{e_i})\cap N_{2k}(\tilde X_{e_{i+1}})\cap \tilde
  X_u^+$ has diameter $\leq \psi'(2k)$.
  Now note that $N_k(F_1)\cap
  N_k(F_2)$ is contained in the $k$-neighborhood of $Q$, hence
  diam$(N_k(F_1)\cap N_k(F_2)) \leq \psi'(2k)+2k$.

   Hence, it suffices to define the function $\phi$ as follows
  $$\phi(k)=\max\{\phi_v(2k)+2k, \psi_e(2k)+2k,  
  \psi'(2k)+2k, \mbox{diam }\tilde X_e+2k\},$$ where the maximum is
  taken over all vertex spaces $\tilde X_v$, all inclusions of edge
  spaces which are lines into G vertex spaces, and all edge spaces
  $\tilde X_e$ which are segments.
\end{proof}

Now let us recall the setting. $L$ is a limit group and it is embedded
in an $\omega$-rft $G=\pi_1(Y)$ of height $n>0$. We view $Y$ as a
graph of spaces resulting from attaching the last stage; that is,
either
\begin{enumerate}
\item $Y=M\sqcup \sqcup_t S^1_t\times [0,1]\sqcup N/\sim$ with $N$ a
  hyperbolic surface, or
\item $Y=M\sqcup S^1\times [0,1]\sqcup N/\sim$ with $N$ a torus. In
this case the attaching circle in $M$ generates a maximal abelian
subgroup.
\end{enumerate}
Let $p:X\to Y$ be the covering
 space with $\pi_1(X)=L$. According to
Theorem \ref{main} $X$ contains a core $C\subset X$. The induced graph
of spaces decomposition of $C$ is a special graph of CAT(0) spaces,
and there are two types of vertex spaces, the $M$-type and the
$N$-type (these are in fact the cores of the covering spaces of $M$
 and $N$). 

We will need several lemmas.

\begin{lemma}\label{max}
Suppose a group $G$ acts geometrically on a CAT(0) space with isolated
flats $(X, \mathcal F, \phi)$ with $\mathcal F$ $G$-invariant. Let $g$
be an element of $G$ that acts as a translation on a periodic geodesic
$\ell \subset X$. If $\ell \parallel F$, for some $F \in \mathcal F$,
then $g\in Stab(F)$ (and hence $g$ does not generate a maximal abelian
subgroup).
\end{lemma}

\begin{proof}
Since $\ell$ is parallel to the flat $F$ there is a neighborhood
$N_k(F)$ that contains $\ell$. It follows that $N_k(gF)$ contains
$g\ell=\ell$, hence $N_k(F)\cap N_k(gF)$ contains $\ell$. Since both
$F$ and $gF$ belong to $\mathcal F$, (2) in the
definition of isolated flats implies $gF=F$, i.e. $g \in Stab(F)$. 
\end{proof}

\begin{lemma}\label{strips}
Suppose that $N$ is a torus.  No two distinct edge spaces incident to
an $M$-type vertex space $V$ in $\tilde C$ cobound a strip in the associated
extended vertex space $V^+$.
\end{lemma}

\begin{proof}
Suppose this is false and let $\ell_1,\ell_2$ be two distinct parallel
lines that are edge spaces incident to $V$. Then $\ell_1,\ell_2$ are
also edge spaces in $\tilde Y$ contained in an extended vertex space
$W^+\supset V^+$. Since the graph of groups associated to $Y$ has one
edge and two vertices, there is $h\in Stab(W)$ with
$h(\ell_1)=\ell_2$.  Let $g_1\in\pi_1(Y)$ be a primitive element with axis
$\ell_1$. Thus $g_2=hg_1h^{-1}$ is a primitive element with axis
$\ell_2$. 
\begin{comment}
There are only finitely many strips with the same width as
$S$ and with one boundary component $\ell_1$ and the other an edge
space. Therefore for some $k>0$ we must have $g_1^k(S)=S$ and in
particular $\ell_2$ is an axis of $g_1^k$. Since $Stab(\ell_2)$ is
cyclic and generated by $g_2$ we must have $g_1^k=g_2^l$ for some
$l$. 
\end{comment}
By Lemma \ref{half-flats0} we conclude that  $g_1^k=g_2^l$ for some
$k$ and $l$. That is to say,
$$g_1^k=hg_1^lh^{-1}$$ 
This is an equation in $Stab(W)$, which is a
limit group. It now follows that $k=l$ and $g_1h=hg_1$ (see
Proposition \ref{basic} (4),(5) -- proof: if
$g_1h\neq hg_1$ let $f:Stab(W)\to F$ be a homomorphism to a free group
that does not kill $[g_1,h]$ and obtain a contradiction in
$F$). However, $h$ is not contained in the cyclic group generated by
$g_1$ since it does not stabilize $\ell_1$. This contradicts the
assumption that $g_1$ generates a maximal abelian group in $Stab(W)$.
\end{proof}

\begin{lemma}\label{half-flats}
Let $\ell\subset \tilde C$ be a periodic line in the preimage of the
core in $\tilde Y$ whose projection to $Y$ is a loop that generates a
maximal abelian subgroup of $\pi_1(Y)$. Then $\ell$ does not bound a
half-flat in $\tilde C$.
\end{lemma}

\begin{proof}
We argue by induction on the height of the tower.

{\bf Case 1.} $\ell$ transversely crosses an edge space. Of course, in
this case $N$ must be a torus or else it is clear that $\ell$ does not
bound a half-flat. We will first argue that $\ell$ bounds
at most two half-flats. Let $S$ be the infinite strip associated with
an edge space that $\ell$ crosses. The intersection $\ell\cap S$ is a
segment and it separates $S$ into two components, say $S_1$ and
$S_2$. Any half-flat bounded by $\ell$ must contain either $S_1$ or
$S_2$. If $P,P'$ are two half-flats bounded by $\ell$ that both
contain say $S_1$ then $P\cap P'$ is a convex subset of $P$ that
contains $\ell\cup S_1$ and thus must equal $P$. This proves that
$P=P'$ and that $\ell$ bounds at most two half-flats. Now suppose that
$\ell$ bounds a half-flat $P$. Let $g$ be an element that acts as a
translation on $\ell$. We now conclude that $g^2(P)=P$. In particular,
$g^{2k}(S)$ are other strips associated to edge spaces and all are
parallel to $S$. It follows that all edge spaces that intersect $P$
are parallel to each other, and two consecutive edge spaces along
$\ell$ that are incident to an $M$-type vertex space contradict
Lemma \ref{strips}.

{\bf Case 2.} $\ell$ is contained in a vertex space $V$. If $V$ is an
$N$-type vertex space then either $N$ is a surface and then $\ell$
cannot bound a half-flat or $N$ is a torus and then the image loop
does not generate a maximal abelian subgroup. So we may assume that
$V$ is an $M$-type vertex space. By induction, $\ell$ does not bound
any half-flats in $V$. It remains to rule out half-flats $P$ with
$\partial P=\ell$ and $P$ intersecting some edge spaces. These
intersections must be lines parallel to $\ell$. Let
$\ell_1,\ell_2,\ldots$ be the (finite or infinite) sequence of lines
of intersection between $P$ and the edge spaces ordered according to
distance from $\ell$. Thus the strip $S$ between $\ell$ and $\ell_1$
is contained in $V$ and the strip between $\ell_1$ and $\ell_2$ (or
the half-flat in case there is no $\ell_2$) is contained in the
adjacent $N$-type vertex space $W$. Thus $N$ is a torus. Let $g_1$ be
a primitive element that translates along $\ell_1$. As in the proof of
Lemma \ref{half-flats0} we see that
$g^k\in <g_1>$. It follows that $g^k$, and therefore\footnote{this is
based on the property of limit groups that if $a,b,c$ are nontrivial
elements and $[a,b]=[b,c]=1$ then $[a,c]=1$, see Proposition
\ref{basic}(3).} $g$, belong to the noncyclic abelian subgroup
$Stab(W)\cong\pi_1(N)$, contradiction.
\end{proof}

\begin{thm}\label{ifp}
  Limit groups act geometrically on CAT(0) spaces with isolated flats
  property.
\end{thm}

\begin{proof}
Let $L$ be a limit group and $L\subset G$ an embedding of $L$ in an
$\omega$-rft $G$. Let $n$ be the height of $G$. If $n=0$ then $L$ is the
free product of free groups, free abelian groups, and hyperbolic
surface groups, and the fact that such groups have isolated flats is
straightforward from the definition. We will now assume $n>0$ and that
the theorem holds for limit groups that can be embedded in
$\omega$-rft's of height $<n$.

Let $p:X\to Y$, $C\subset X$ be as in our setting in this section. We
will verify the conditions of Theorem \ref{combination} for the graph
of spaces decomposition of $C$ inherited from $X$. If $N$
is a hyperbolic surface then the $N$-type vertices cannot contain
flats, half-flats, or strips, so declare that N-type vertices are G
vertices and M-type vertices are B vertices. Now assume that $N$ is a
torus. Then declare that M-type vertices are G vertices and N-type
vertices are B vertices. That $M$-type vertices satisfy (1) and (2)
follows from Lemmas \ref{strips} and \ref{half-flats}.
\end{proof}

\section{Appendix: Sketch of proof of Theorem \ref{embedding}}

In this Appendix we outline a proof of Theorem \ref{embedding} due to
Sela. The reader is assumed to have some familiarity with limit groups
e.g. as in the expository papers
\cite{champetier-guirardel},\cite{paulin},\cite{bf:notes} as well as
with the language of JSJ decompositions.

First, there is a basic fact about $\omega$-rft's easily proved by
induction on height.

\begin{lemma}\label{maxab}
Let $Y$ be an $\omega$-rft and $A$ a noncyclic abelian subgroup
of $\pi_1(Y)$. Then $A$ is conjugate into the fundamental group of a
torus added at some stage in the construction of $Y$.
\end{lemma}

It then follows that in the definition of an $\omega$-rft we could add
another building block
\begin{enumerate}
\item[(T)] (Torus block) $Y_n=Y_{n-1}\sqcup T^k\times [0,1]\sqcup
  T^{l}/\sim$ where $T^k$ and $T^l$ are $k$- and $l$-tori, $k<l$,
  $T^k\times\{1\}$ is identified with a coordinate torus in $T^l$, and
  $T^k\times\{0\}$ is glued to $Y_{n-1}$ by a map that is
  $\pi_1$-injective and its image is a maximal abelian subgroup.
\end{enumerate}
If $k=1$ this is exactly building block (A) and if $k>1$ then it
follows from Lemma \ref{maxab} that the gluing map to $Y_{n-1}$ can be
taken to be an isomorphism to a previously added torus. Thus we could
modify this earlier step and glue in the $l$-torus instead of the
$k$-torus and the modified $\omega$-rft $Y_{n-1}'$ would be homotopy
equivalent to $Y_n$.

Let $L$ be a given limit group. We want to show that $L$ embeds in
some $\omega$-rft. 
 
An important fact about limit groups is that any sequence of
epimorphisms $L_1\to L_2\to\cdots$ eventually consists of
isomorphisms\footnote{Proof: The sequence of algebraic varieties
$Hom(L_1,SL_2(\R))\supset Hom(L_2,SL_2(\R))\supset\cdots$ eventually
consists of equalities.}. We can thus use induction and assume that all
proper limit group quotients of $L$ embed in $\omega$-rft's.

If $L=L_1*L_2$ is a free product and $L_i\subset \Gamma_i$ ($i=1,2$)
are embeddings in $\omega$-rft's then $L=L_1*L_2\subset
\Gamma_1*\Gamma_2$ and the free product of $\omega$-rft's is an
$\omega$-rft. Thus we can assume that $L$ is freely indecomposable.

If $L$ is abelian it is already an $\omega$-rft. Otherwise, one
carefully constructs a proper quotient $q:L\to L'$ with certain
properties. This amounts to finding a ``strict quotient'' (see
\cite[5.9,5.10]{zlil1}). We state this result as it appears in
\cite[1.14,1.25]{bf:notes}.

\def\fg{finitely generated}
\def\limitgroup{L}
\def\gad{generalized abelian decomposition}
\def\QH{QH}
\def\blackbox{V}
\begin{prop}\label{CLG}
The class of limit groups coincides with the following hierarchy of groups
(Constructible Limit Groups)

Level 0 of the hierarchy consists of \fg\ free groups, \fg\ free
abelian groups, and fundamental groups of closed surfaces that support
pseudoAnosov homeomorphisms.

A group $\limitgroup$ belongs to level $\leq n+1$ iff either it has a
free product decomposition $\limitgroup=\limitgroup_1*\limitgroup_2$
with $\limitgroup_1$ and $\limitgroup_2$ of level $\leq n$ or it has a
homomorphism $\rho:\limitgroup\to\limitgroup'$ with $\limitgroup'$ of
level $\leq n$ and it has a \gad\ such that
\begin{itemize}
\item $\rho$ is injective on the peripheral subgroup of each abelian
vertex group.
\item $\rho$ is injective on each edge group $E$ and at least one of
  the images of $E$ in a vertex group of the
  one-edged splitting induced by $E$ is a maximal abelian subgroup.
\item The image of each $\QH$-vertex group is a non-abelian subgroup of
$L'$.
\item For every rigid vertex group $V$, $\rho$
is injective on the ``envelope'' $\tilde\blackbox$ of $\blackbox$,
defined by first replacing each abelian vertex with the peripheral
subgroup and then letting $\tilde\blackbox$ be the subgroup of the
resulting group generated by $\blackbox$ and by the centralizers of
incident edge-groups.
\end{itemize}
\end{prop}

\begin{remark}
Sela's proof that every limit group belongs to this hierarchy goes
like this. We may assume $L$ is freely indecomposable. Consider the
abelian JSJ decomposition of $L$. Let $f_i:L\to F$ be a sequence of
homomorphisms to a free group so that each $1\neq x\in L$ is killed by
only finitely many $f_i$'s. Let $g_i$ be obtained from $f_i$ by
conjugating and composing with the elements of the modular group of
the JSJ decomposition so that the sum of the lengths of the images of
fixed generators of $L$ is as small as possible. A subsequence of the
$g_i$'s converges to an action of $L$ on an $\R$-tree and $q:L\to L'$
divides out the kernel of the action. One then checks the properties
stated in Proposition \ref{CLG}.
\end{remark}

So now suppose that we have a freely indecomposable, nonabelian limit
group $L$. Let $q:L\to L'$ be as above and let $i:L'\hookrightarrow
\Gamma'$ be an embedding in an $\omega$-rft $\Gamma'$. The idea is to
construct an $\omega$-rft $\Gamma$ by attaching things to $\Gamma'$ so
that $L$ embeds in $\Gamma$. This is explained in
\cite[1.12]{zlil2}. Rather than repeat the definition we illustrate
the construction on examples that cover the cases when the JSJ
decomposition has one edge. The reader can easily extrapolate the
general case (or refer to \cite[1.12]{zlil2}).

\begin{example}
Suppose $L=A*_E B$ with $A,B$ rigid. We have the composition
$\nu=iq:L\to \Gamma'$ 
which embeds $A$ and $B$ and their envelopes. Let $U$ be the maximal
abelian subgroup of $\Gamma'$ that contains $\nu(E)$. Form $\Gamma$ as 
$$\Gamma=\Gamma'*_U (U\times \Z)$$ This corresponds to step (T) of the
construction of $\omega$-rft's. Let $t$ be a generator of $\Z$. Now
define $$j:L\to \Gamma$$ by $j|A=\nu|A$ and $j|B=(\nu|B)^t$
(i.e. $j(b)=t\nu(b)t^{-1}$ for $b\in B$). 

It is interesting to note that a special case of this construction was
discovered by G. Baumslag \cite{MR25:3980} in 1962!
\end{example}

\begin{example}
Let $L=A*_E$ with $A$ rigid. We again have $\nu:L\to \Gamma'$ and form
$\Gamma$ similarly
$$\Gamma=\Gamma'*_U (U\times\Z)$$
where $U$ is the maximal abelian subgroup of $\Gamma'$ that contains
$\nu(E)$. Then define $j:L\to \Gamma$ by $j|A=\nu|A$ and $j$ sends the
``stable letter'' $s$ to $t\nu(s)$ where $t$ is a generator of $\Z$.
\end{example}

\begin{example}
Let $L=A*_E B$ with $B$ rigid and $A$ abelian. Then perform step (T)
to the $\omega$-rft $\Gamma'$ by gluing an abelian group to the maximal
abelian subgroup containing the image of $E$, and send $A$ there.
\end{example}

\begin{example}
Let $L=A*_EB$ where $B$ is rigid and $A$ is a QH-vertex group and $E$
corresponds to the boundary of the surface. Let $\nu:L\to \Gamma'$ be the
composition $L\to L'\to \Gamma'$ as above and let $E'=\nu(E)$. Now apply
step (Q) and form
$$\Gamma=\Gamma'*_{E'}A'$$
where $(A',E')$ is a copy of $(A,E)$. Define the embedding $j:L\to \Gamma$ by
$j|B=\nu|B$ and $j:A\to A'$ is the identifying isomorphism.
\end{example}

\begin{remark}
We outlined a construction of an embedding $j:L\to\Gamma$ of a given
limit group $L$ into an $\omega$-rft $\Gamma$. This embedding depends
on the choice of a sequence of strict quotients defined on the free
factors of $L$ and the free factors of subsequent quotients ending in
free groups. Such a sequence is called a {\it strict resolution} and
$j:L\to\Gamma$ is the associated {\it completion}. Every limit group
has a {\it canonical} finite collection of resolutions and each
resolution gives rise to a completion $L\to\Gamma$ (which may not be
an embedding if the resolution is not strict). Thus every limit group
has a canonical finite collection of completions. The study of
completions is the main topic of \cite{zlil2}.
\end{remark}

\bibliographystyle{amsalpha}

\bibliography{ref}

\end{document}